\documentclass[12pt]{article}

\usepackage{amstext,amsmath,amssymb,amsthm}
\usepackage{graphicx,subcaption}

\theoremstyle{theorem}

\newtheorem{proposition}{Proposition}
\newtheorem{lemma}{Lemma}
\theoremstyle{definition}

\newtheorem*{remark}{Remark}

\newcommand{\vs}{\vspace{0.5cm}}
\newcommand{\R}{\mathbb{R}}

\begin{document}

\title{Revisiting second-order optimality conditions 
for equality-constrained minimization problem}
\author{Luca Amodei}


\maketitle

\begin{abstract} The aim of this note is to give a geometric insight into  the classical second order optimality conditions for equality-constrained minimization problem. We show that the 
Hessian's positivity of the Lagrangian function associated to the problem  at a local minimum point $x^*$ corresponds to inequalities between the respective algebraic curvatures at point $x^*$ 
of the  hypersurface $\mathcal{M}_{f, x^*}=\{ x \in \R^n \, | \, f(x) = f(x^*)\}$ defined  by the objective function $f$   and the submanifold $\mathcal{M}_g = \{ x \in \R^n \, | \, g(x)= 0 \}$ defining the contraints. These inequalities highlight a geometric evidence on how, in order  to guarantee the optimality,  the submanifold $\mathcal{M}_g$ has to be locally included in the half space  $\mathcal{M}_{f, x^*}^+ = \{ x \in \R^n \, | \, f(x) \geq f(x^*)\}$ limited by the hypersurface 
$\mathcal{M}_{f, x^*}.$  This presentation can be used for educational purposes and help to a better understanding of this property.
\end{abstract}

\section{The classical second-order conditions}

We consider the optimization problem with equality constraints : 

\begin{equation} \begin{array}{cc} \min & f(x) \\ x \in \R^n & \\
\textsl{ and }  g(x) = 0 & \end{array} \label{eq1}
\end{equation}

defined by the objective function $f$ and  the vector function $g =(g_1, \cdots, g_m), \, m \leq n,$ for the constraints. We suppose that $f$ and $g$ are $C^2$ differentiable functions defined on $\R^n.$ \footnote{In the same way, we could consider that the functions $f$ and $g$ are defined on an open set $\Omega \subset \R^n, \, \Omega \neq \emptyset.$}
\vs
  
We suppose that $0 \in \R^m$ is a regular value of  $g$, i.e. the Jacobian 
$J g(x)$ is full rank for every $x \in g^{-1}(\{0\}).$ 
Thus, $\mathcal{M}_g= g^{-1}(\{0\})= \{ x \in \R^n \, | \, g(x)= 0 \}$ is a closed embedded $(n-m)$-dimensional submanifold of $\R^n,$ and we know that, for every $x \in \mathcal{M}_g,$ the tangent space $T_{x} \mathcal{M}_g$ at  $x$ to $\mathcal{M}_g$   is equal to  
$\textrm{Ker}(J g(x)).$

\vs 

We recall the classical first-order and second-order optimality  conditions for the problem \eqref{eq1} (see ~\cite{LY} and ~\cite{NW}).

If the point $x^* \in \mathcal{M}_g$ is a local minimum of \eqref{eq1}, the first-order necessary conditions are  given by  
\begin{equation}  \nabla f(x^*) - \sum_{i=1}^m \lambda_i^* \nabla g_i(x^*) = \nabla f(x^*) - J g(x^*)^T \lambda^* =  0, \label{eq2} \end{equation}
where $\nabla f(x^*), \nabla g_i(x^*), i=1, \ldots , m,$ are  the gradients of $f$ and $g_i$ at  $x^*$,  $J g(x^*)$ the Jacobian matrix of $g$ at $x^*,$ and $\lambda^*=(\lambda_1^*, \cdots, \lambda_m^*) \in  \R^m$ the  Lagrange multipliers vector.  

Equation \eqref{eq2} says that the orthogonal of the tangent spaces 
$T_{x^*} \mathcal{M}_{f,x^*}$ and $T_{x^*} \mathcal{M}_g$ at $x^*$ satisfy the inclusion 
$$(T_{x^*} \mathcal{M}_{f,x^*})^\perp \subseteq (T_{x^*} \mathcal{M}_g)^\perp,$$   
or equivalently 
$$T_{x^*} \mathcal{M}_g \subseteq  T_{x^*} \mathcal{M}_{f,x^*}.$$   

Notice that for $m=1,$ the tangent spaces at $x^*$ are equal.

By introducing the Lagrangian function $L(x,\lambda)$ associated to the constrained optimization problem \eqref{eq1} 
$$ L(x,\lambda) = f(x) - \sum_{i=1}^m \lambda_i g_i(x),$$ the conditions \eqref{eq2} are given by   
$$\nabla_x L(x^*, \lambda^*) = 0,$$ where $\nabla_x$ is the gradient of $L$ with respect to $x.$

The second-order necessary conditions are the following (see ~\cite{LY} and ~\cite{NW}) :  
\begin{equation} v^T \left(\nabla^2 f(x^*) - \sum_{i=1}^m \lambda_i^* \nabla^2 g_i(x^*) \right) v \geq 0, \, \forall 
v \in \textrm{Ker}(J g(x^*)),
\label{eq3} 
\end{equation}
with $\nabla^2 f (x^*) $ and $\nabla^2 g_i (x^*) $ the Hessian matrices of  $f$ and $g_i$ at $x^*,$ and supposing that first-order conditions \eqref{eq2} are verified. 

\vs 

Again, the Lagrangian function $L$ allows to formulate concisely conditions \eqref{eq3} : 
$$v^T \nabla^2_{xx} L(x^*, \lambda^*) v  \geq 0, \, \forall 
v \in \textrm{Ker}(J g(x^*)),$$ where $\nabla^2_{xx} L$ is the Hessian of $L$ with respect to $x.$    

Moreover, conditions \eqref{eq3} with strict inequality for all  $v \in \textrm{Ker}(J g(x^*)), v \neq 0,$ along with first-order conditions  \eqref{eq2},  are sufficient conditions  for $x^* \in \mathcal{M}_g$ to be a local minimum of the problem \eqref{eq1}.   

\vs

In an introductory course in optimization it makes sense to illustrate the first-order conditions \eqref{eq2} with a drawing, at least  in the simple case $ n = 2 $ and $ m = 1 $ (curves on the plane): the two curves $\mathcal{M}_{f, x^*} $ and $\mathcal{M}_g $ in the plane $\R^2$ have a common tangent at point $ x ^ *. $  
Naturally follows the question for second-order conditions : how conditions \eqref{eq3}  impose  a new geometric contraint to first-order equation \eqref{eq2} ? How should the picture be in the simple case $ n = 2 $ and $ m = 1 $ ? Furthermore, since second-order conditions \eqref{eq3} with strict inequality combined with first-order conditions \eqref{eq2} also give sufficient conditions for $x^* \in \mathcal{M}_g$ to be a local minimum, it must garantee that $\mathcal{M}_g$ is locally on the half space  $\mathcal{M}_{f, x^*}^+ = \{ x \in \R^n \, | \, f(x) \geq f(x^*)\}.$

Let us start by considering  the second-order conditions \eqref{eq3} in the case $n=2$ and $m=1 :$ 
\begin{equation} v^T \nabla^2 f(x^*) v  \geq \lambda^* \, v^T \nabla^2 g(x^*) v, \label{eq33} \end{equation} for every tangent vector 
$v \in \R^2$ of the curve $\mathcal{M}_g$ at point $x^*.$ We can take $v=v_g = (\partial_{x_2} g(x^*), -\partial_{x_1} g(x^*))^T$ or $v= v_f =     
(\partial_{x_2} f(x^*), -\partial_{x_1} f(x^*))^T= \lambda^* v_g,$ since first-order conditions \eqref{eq2} are satisfied. We suppose that $\nabla f(x^*) \neq 0, $ and consider the unit vectors 
$u_g = v_g / \| v_g \|$ and $u_f = v_f / \| v_f \|.$  

A classical result of  elementary differential geometry (see ~\cite{G}) gives the algebraic curvatures $\kappa_f(x^*)$ and $\kappa_g(x^*)$ of the curves $\mathcal{M}_{f, x^*} $ and $\mathcal{M}_g $ at a regular point $x^*$ :
$$\kappa_f(x^*) = - \frac{u_f^T \nabla^2 f(x^*) u_f}{\| \nabla f(x^*) \|} \, \textrm{ and }  \kappa_g(x^*) = - \frac{u_g^T \nabla^2 g(x^*) u_g}{\| \nabla g(x^*) \|}.$$    
Now, if we divide inequality \eqref{eq33} by $\| \nabla f(x^*) \|$,  take $v = u_f$ and use first-order conditions \eqref{eq2}, we obtain 
\begin{equation} - \frac{u_f^T \nabla^2 f(x^*) u_f}{\| \nabla f(x^*) \|}  \leq - 
\lambda^*  \frac{u_f^T \nabla^2 g(x^*) u_f}{\| \nabla f(x^*) \|} = 
- \frac{\lambda^*}{| \lambda^* |} \,   
\frac{u_g^T \nabla^2 g(x^*) u_g}{\| \nabla g(x^*) \|}.  \label{eq34} \end{equation}
Therefore, the algebraic curvatures satisfy 
\begin{equation} \kappa_f(x^*) \leq \pm \kappa_g(x^*), \label{eq35} \end{equation} 
with $+1$ if the vectors $\nabla f(x^*)$ and  $\nabla g(x^*)$ have the same orientation and $-1$ otherwise.  
We can easily verify that inequality \eqref{eq35} on the curvatures  is necessarily  satisfied if  $\mathcal{M}_g$ is locally on the half space  $\mathcal{M}_{f, x^*}^+.$  Figure 1 gives, according to the sign of the algebraic curvatures, the four situations for which the curves $\mathcal{M}_{f, x^*} $ and $\mathcal{M}_g $ satisfy inequality \eqref{eq35}.     

\begin{figure}
\centering
\begin{subfigure}[b]{0.4\linewidth}
    \includegraphics[width=\linewidth]{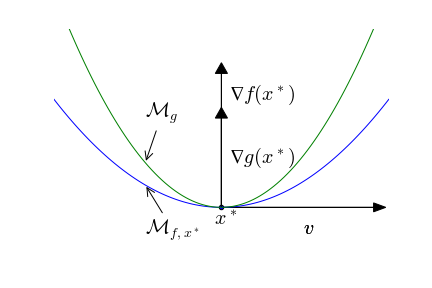}
    \caption{$\kappa_f(x^*) \geq 0$ and $\kappa_g(x^*) \geq 0.$}
  \end{subfigure}
  \begin{subfigure}[b]{0.4\linewidth}
    \includegraphics[width=\linewidth]{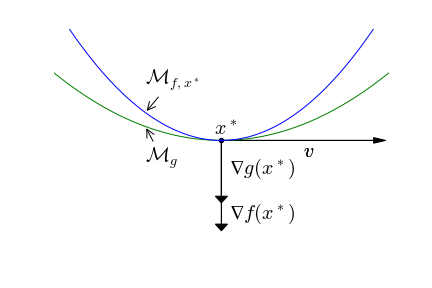}
    \caption{$\kappa_f(x^*) \leq 0$ and $\kappa_g(x^*) \leq 0.$}
  \end{subfigure}
  \begin{subfigure}[b]{0.4\linewidth}
    \includegraphics[width=\linewidth]{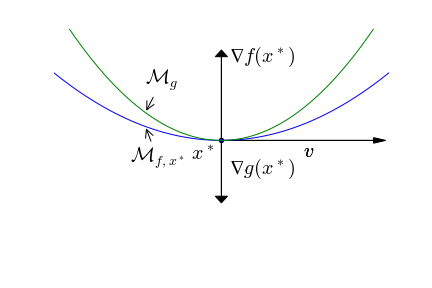}
    \caption{$\kappa_f(x^*) \geq 0$ and $\kappa_g(x^*) \leq 0.$}
  \end{subfigure}
  \begin{subfigure}[b]{0.4\linewidth}
    \includegraphics[width=\linewidth]{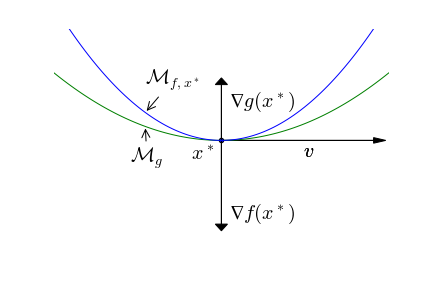}
    \caption{$\kappa_f(x^*) \leq 0$ and $\kappa_g(x^*) \geq 0.$}
  \end{subfigure}
  \caption{Algebraic curvatures of $\mathcal{M}_{f,x^*}$ and $\mathcal{M}_g.$}
\end{figure}

To my knowledge, in the classical litterature on the subject,  inequality \eqref{eq35} is not explicitly presented despite its intuitive and natural aspect.   

In fact inequality \eqref{eq35} and more generally inequality \eqref{eq3}  directly appears  by requiring  positivity of the  second derivative of the function $ (f \circ \gamma)(t)$ at $t=0$ (like for unconstrained minimization), for every regular curve $\gamma(t)$ in $\mathcal{M}_g$  such that $\gamma(0) = x^*$ (see ~\cite{LY}). 

By elementary calculus we obtain  
$$ (f\circ \gamma)''(0) = \gamma'(0)^T \, \nabla^2 f(x^*) \, \gamma'(0) + \nabla f(x^*)^T \, \gamma''(0).$$ 
In the particular case $m=1$, the constraint  $\gamma(t) \in \mathcal{M}_g,  \, \forall t,$ gives    
\begin{equation} (g \circ \gamma)''(0) = \gamma'(0)^T \nabla^2 g (x^*) \gamma'(0) + 
\nabla g (x^*)^T \gamma''(0) = 0.\label{eq41} \end{equation}    

Combining inequality  $(f \circ \gamma)''(0) \geq 0$ with \eqref{eq41}, dividing  
by $\| \nabla f(x^*) \|,$ and using first-order conditions \eqref{eq2}, we obtain  

\begin{equation} - \, \frac{\gamma'(0)^T \, \nabla^2 f(x^*) \, \gamma'(0)}{\| \nabla f(x^*) \|} \, \leq  
\, - \frac{\lambda^*}{| \lambda^* |} \, \frac{\gamma'(0)^T \nabla^2 g (x^*) \gamma'(0)}{\| \nabla g(x^*) \|} \label{eq43},\end{equation} 
which is  equivalent to \eqref{eq35} in the particular case $n=2.$  
 
The purpose of this note is to show how \eqref{eq43} can be seen, in the case $m=1$,  as an inequality on the respective curvatures of the intersections of $\mathcal{M}_{f, x^*}$ and $\mathcal{M}_g$  with the affine plane 
$\mathcal{L} = \{x^* +  a \gamma'(0)  + b \nabla f(x^*), \, \forall  (a,b) \in \R^2 \}, $ and extend appropriately this property for any value of the codimension $m, \, 
1 \leq m \leq n.$

\section{The curvature of the submanifolds $\mathcal{M}_\MakeLowercase{f, x^*}$ and $\mathcal{M}_\MakeLowercase{g}$}

The curvature of a manifold is a fundamental concept in differential geometry. It describes 
 the way a curve or a surface deviates from being a straight line or a flat plane. A distinction exists between the notion of extrinsic curvature  adapted for submanifolds embedded in an euclidean space and the notion of intrinsic curvature defined for Riemannian manifolds without reference to an ambiant euclidean space. 
As  we deal with the submanifolds $\mathcal{M}_{f, x^*}$ and $\mathcal{M}_g$, we will consider the first notion of curvature.

In this paragraph we introduce  the notion of second fundamental form and show its essential relation to the curvature (see ~\cite{RS22}) of a submanifold. Then, we determine the second fundamental forms of the submanifolds $\mathcal{M}_{f, x^*}$ and $\mathcal{M}_g$ at point $x^*.$ Finally we show how second-order necessary conditions \eqref{eq3} correspond to inequalities between the curvatures  of the submanifolds $\mathcal{M}_{f, x^*}$ and $\mathcal{M}_g$ at point $x^*.$    

\subsection{The second fundamental form}

Let $\mathcal{M} \subset \R^n$ be a submanifold. 
For $x \in \mathcal{M},$ there are different ways to define the  second fundamental form at $x.$ For instance, when $M$ is an hypersurface, this can be done by using the Weingarten map (also called the shape operator - see ~\cite{T}). We will use a more general presentation based on the orthogonal projection 
$\Pi(x)$ onto the tangent space $T_x \mathcal{M}$ (see ~\cite{RS22}).
This formulation is particularly adapted to compute the second fundamental form on the submanifolds $\mathcal{M}_{f,x^*}$ and $\mathcal{M}_g$ which are defined implicitly.    

\vs

Let us consider the orthogonal projection $\Pi(x)$ onto the tangent space $T_x \mathcal{M}.$ 
For each $x \in \mathcal{M},$ $\Pi(x) \in \R^{n \times n}$  verifies 
$$\Pi(x)=\Pi(x)^2 = \Pi(x)^T,$$ and the  equivalence 
$$\Pi(x) v = v \Leftrightarrow v \in T_x \mathcal{M}.$$

The second differential form is given by differentiating the map $\Pi : \mathcal{M} \rightarrow \R^{n \times n}$ at point $x.$ 

The derivative $$ d\Pi(x) : T_x \mathcal{M} \rightarrow  \R^{n \times n},$$ 
of $\Pi$ at a point $x$  
 is defined by 
$d\Pi(x)v := (\Pi \circ \gamma)'(0),$ where $\gamma : I=]-\epsilon, \epsilon[ \rightarrow \R^n,$ (with $\epsilon >0 $)  is a parametrized curve such that $\gamma(I) \subset \mathcal{M},$ $\gamma(0) = x,$ and $\gamma'(0) = v.$ 
  
$d\Pi(x) v$ is a matrix and can therefore be multiplied by a vector $u \in \R^n.$ We easily check  the map $(d\Pi(x)v)u $  satisfies the following 
properties. 

\begin{proposition}

For all $x \in \mathcal{M}$ and $u,v \in T_x \mathcal{M},$ we have  
\begin{enumerate} 
\item     
$$(d\Pi(x)v)u \in (T_x \mathcal{M})^\bot, $$
\item 
$$(d\Pi(x)v)u = (d\Pi(x)u)v.$$
\end{enumerate}

\end{proposition}

\begin{proof} See ~\cite{RS22}, Chapter 3. 
\end{proof}

The family of  symmetric bilinear maps 
$$h_x : T_x \mathcal{M} \times T_x\mathcal{M} \longrightarrow T_x \mathcal{M}^\perp,$$   
defined for every $x \in \mathcal{M},$ by $h_x(u,v) = (d\Pi(x)u)v,$ with $u,v \in T_x \mathcal{M},$   is called the second fundamental form on $\mathcal{M}.$

\vs

\begin{remark}  The terminology ``second fundamental form ''  suggests that the value $h_x(u,v)$ is a scalar. In fact, to be consistent with the terminology,  
the second fundamental form is generally defined as the scalar product $\left\langle  h_x(u,v), n \right\rangle$ of $h_x(u,v)$ with a vector $n \in T_x \mathcal{M}^\perp.$ $\left\langle  h_x(u,v), n \right\rangle$  is then called second fundamental form along the normal vector $n$ (see ~\cite{DC}). 
\end{remark}

\vs 
The following result gives the fundamental relation between the second fundamental form $h_x$ and the curvature at a point $x \in \mathcal{M}.$  

\begin{proposition}

Consider $x \in \mathcal{M},$ a unit tangent vector 
$v \in T_x \mathcal{M}, $ and $\mathcal{L}$ the affine space defined by  
$$\mathcal{L}_v = \{ x + tv + w \, | \, t \in \R, \, w \in T_x \mathcal{M}^\perp \}.$$
Let  $\gamma$   
$$\gamma: ]\epsilon, \epsilon[ \rightarrow \mathcal{M} \cap \mathcal{L}_v,$$ 
a second-order differentiable arc length parametrized curve
 which verifies  $\gamma(0)=x, \, \gamma'(0) = v, $ and $\| \gamma'(s) \| = 1, \, \forall s \in 
]\epsilon, \epsilon[,$ ($\epsilon > 0$).

Then $$\gamma ''(0) = h_x(v,v).$$
\end{proposition}

\vs 

\begin{proof} See ~\cite{RS22}, Chapter 3. \end{proof}

The existence of such curves $\gamma$ will be given in section 3.  

\vs 
\subsection{Inequality of the curvatures} 

Let us determine $h_{x^*}(v,v)$  for $\mathcal{M} = \mathcal{M}_{f,x^*}=\{ x \in \R^n \, | \, f(x) = f(x^*)\}$   and $\mathcal{M} = \mathcal{M}_g=\{ x \in \R^n \, | \, g(x) = 0 \}.$ 
To simplify the notations we consider a generic point $x \in \mathcal{M}$ and $\mathcal{M} = \mathcal{M}_f=\{ x \in \R^n \, | \, f(x) = 0\}$ instead of  $\mathcal{M}_{f,x^*}.$ This last choice is supported by the fact that the second fundamental form, as we are going to see, is identical for the hypersurface $\mathcal{M}_f$ definied by           
a function $f$ and the hypersurface  $\mathcal{M}_{f,x^*}= \mathcal{M}_{f - f(x^*)}$ defined by the translated function $f - f(x^*).$

Let us consider $\Pi_f(x)$ the orthogonal projection onto $T_x \mathcal{M}_f.$ We suppose that  $x$ is a regular point.   
The  vector $\nabla f(x)$ generates the normal space $(T_x \mathcal{M}_f)^\perp.$  By taking the normalized vector $\nu_f(x) = \frac{\nabla f(x)}{\| \nabla f(x) \|}$, the orthogonal projection  
onto $(T_x \mathcal{M}_f)^\perp$ is given by the matrix $\nu_f(x) \, \nu_f(x)^T$ where $\nu_f(x)$ is a column vector.  By complementarity we obtain : 
$$\Pi_f(x) = I_n - \nu_f(x) \, \nu_f(x)^T.$$  

For the orthogonal projection $\Pi_g(x)$ onto $T_x \mathcal{M}_g,$ we have to consider the transposed Jacobian matrix $J g(x)^T$ whose columns generate the normal space $(T_x \mathcal{M}_g)^\perp$ at point $x.$ The orthogonal projection onto   $(T_{x} \mathcal{M}_g)^\perp$ is then given by 
$$J g(x)^T \, (J g(x) \, J g(x)^T)^{-1} \, J g(x),$$ from which we deduce  
$$ \Pi_g(x) = I_n - J g(x)^T \, (J g(x) \, J g(x)^T)^{-1} \, J g(x).$$  
         
Let us determine the second fundamental forms. 

By differentiating $\nu_f(x)$, taking $v \in T_x \mathcal{M}_f$ and using $\nu_f(x)^T v = 0,$ we get 
$$(d \Pi_f(x)v)v = h_{f,x}(v,v) = - \nu_f(x) \frac{v^T \nabla^2 f(x)v} {\| \nabla f(x)\| }.$$
The second fundamental form at $x \in \mathcal{M}_f$ along the normal vector $\nu_f(x)$ is thus given by \begin{equation} \langle \nu_f(x), h_{f,x}(v,v) \rangle = \nu_f(x)^T  h_{f,x}(v,v) = -  \frac{v^T \nabla^2 f(x)v} {\| \nabla f(x)\| }. \label{eq6} \end{equation}

In the same way, by differentiating $J g(x),$ taking $v \in T_x \mathcal{M}_g$ and using $J g(x) v = 0,$ we get 
$$(d \Pi_g(x)v)v = h_{g,x} (v,v) = - J g(x)^ T \, 
 (J g(x) J g (x)^T)^{-1} \, v^T \nabla^2 g(x) \, v,$$ where 
$v^T \nabla^2 g(x)  \, v $ is the $\R^m$ column vector  with coefficients $v^T \nabla^2 g_i(x)  \, v, \, i=1, \ldots, m.$   

\vs

The second fundamental form at $x \in \mathcal{M}_g$ along the normal vector $\nu_f(x)$ is then  given by 
\begin{equation} 
\begin{split} 
\langle \nu_f(x), h_{g,x}(v,v) \rangle  &= \nu_f(x)^T  h_{g,x}(v,v) 
\\  &= - \nu_f(x)^T J g(x)^ T \, 
 (J g(x) J g (x)^T)^{-1} \, v^T \nabla^2 g(x) \, v.  
\end{split} \label{eq4}
\end{equation}

We can state the main result. 

\vs 
\begin{proposition}

If $x^* \in \mathcal{M}_g$ is a local minimum, the second-order optimality conditions \eqref{eq3} are equivalent to 
\begin{equation}  \langle \nu_f(x^*), h_{f,x^*}(v,v) \rangle \leq \langle \nu_f(x^*), h_{g,x^*}(v,v) \rangle, \;  \forall v \in \textrm{Ker}(J g(x^*)) = T_{x^*} \mathcal{M}_g. \label{eq7} \end{equation} 
  
\end{proposition}
\begin{proof}  

By multiplying equality \eqref{eq2} on the left by $J g(x^*),$ we get 
$$J g(x^*)\, \nabla f(x^*) - (J g(x^*) J g(x^*)^T) \lambda^* =  0,$$ and thus 
\begin{equation} \lambda^* = (J g(x^*) J g(x^*)^T)^{-1} \,  J g(x^*) \, \nabla f(x^*). \label{eq55} \end{equation}   
Equality \eqref{eq4} becomes  therefore
\begin{equation} \begin{split} \langle \nu_f(x^*), h_{g,x^*}(v,v) \rangle &= 
- \frac{1}{\| \nabla f(x^*)\| } 
(\lambda^*)^T \, v^T \nabla^2 g(x^*) \, v \\ &=  - \frac{1}{\| \nabla f(x^*)\| } 
 v^T \left( \sum_{i=1}^m \lambda^*_i \nabla^2 g_i(x^*) \right) v. \end{split} \label{eq5} \end{equation}
We otain the result by dividing inequality \eqref{eq3} by $\| \nabla f(x^*) \|$ and using equalities \eqref{eq6} and \eqref{eq5}.    
 
\end{proof}

\begin{remark} 

\begin{enumerate}

\item 
Inequality \eqref{eq7} is the generalization of \eqref{eq35} to any values of the dimension $n$ of the space and the codimension $m$ ($1 \leq m \leq n)$  of the submanifold $\mathcal{M}_g.$  

\item 
Since second-order conditions assume first-order conditions, we have $v \in T_{x^*} \mathcal{M}_{f,x^*}$ for any  $v \in T_{x^*} \mathcal{M}_g.$  

\item 
Proposition 3.2.1 in ~\cite{R} also presents the second-order optimality conditions using the second fundamental forms. However, the  parametric form used in ~\cite{R} to represent the submanifolds seems inappropriate.
Moreover, the presentation does not highlight the geometric interpretation in terms of curvature inequalities.
\end{enumerate}
\end{remark}

\vs

In order to complete the discussion we have to show the local existence of curves $\gamma$  in the submanifolds $\mathcal{M}_{f, x^*}$ and $\mathcal{M}_g$  verifying the hypotheses of proposition 2. 
This is obtained in the following paragraph by using implicit function theorem.  Moreover, in section 4 a  geometric proof of the theorem on second-order sufficient optimality conditions is given.

 \section{Implicit curves on $\mathcal{M}_\MakeLowercase{f, x^*}$ and $\mathcal{M}_\MakeLowercase{g}.$}

We consider first the case of the hypersurface  
$\mathcal{M}_{f, x^*}=\{ x \in \R^n \, | \, f(x) = f(x^*)\}.$ We suppose that $x^* \in \R^n$ is a regular point of $f,$ i.e. $\nabla f(x^*) \neq 0.$ 

Let $V \in \R^{n \times (n-1)}$ a matrix whose columns $v_1, \cdots, v_{n-1}$ generate the tangent space $T_{x^*} \mathcal{M}_{f, x^*}.$   
We define the function $\tilde{f} : \R^{n-1} \times\R  \rightarrow \R$ by 
$$\tilde{f}(a,b) = f(x^* + Va +  b \nabla f(x^*) ),$$
for  $a \in \R^{n-1}$ and $b \in \R.$  
Since $\partial_b \tilde{f}(0,0) = \nabla f(x^*)^T \nabla f(x^*) \neq 0,$ by using the implicit function theorem, we can define an open ball  $B(0, r) \subset \R^{n-1}$ of  radius $r > 0,$  
a centered open interval $I_{r'}= ]-r',r'[ \subset \R$  (with $r' > 0$),  and an implicit function $\varphi : B(0,r) \rightarrow I_{r'}$ such that 
$$(a,b) \in B(0,r) \times I_{r'} \; \textrm{ and } \; \tilde{f}(a,b) = \tilde{f}(0,0) = f(x^*)  \Leftrightarrow a \in B(0,r) \; \textrm{ and } \; b= \varphi(a).$$ 
The implicit function $\varphi$ satisfies $\varphi(0) = 0,$ and $\nabla \varphi(0) = 0.$ 
The first equality is due to the unicity of the implicit function. The second is deduced from  the 
equality  $ \nabla \varphi(0)^T = - 
(\nabla f(x^*)^T \nabla f(x^*))^{-1} \, \nabla f(x^*)^T V ,$ satisfied by the gradient of the implicit function $\varphi$  and $\nabla f(x^*)^T V = 0,$  since the  columns of $V$ are  vectors of  $T_{x^*} \mathcal{M}_{f, x^*}.$ 
 
The restriction of the implicit function $\varphi$ to a line $\R v$ generated by a fixed tangent vector $v \in T_{x^*} \mathcal{M}_{f, x^*}$ defines an implicit curve of $\mathcal{M}_{f, x^*}.$ More precisely, if $v = Va$ for some $a \in \R^{n-1}, \, a \neq 0,$ we can define the function $\varphi_a(t) = \varphi(ta)$ for all $t \in \R$ such that $ta \in B(0,r).$ 
Let $\mathcal{L}_{f,v}$ the normal section defined by the tangent vector $v,$ 
$\mathcal{L}_{f,v} = \{ x^* + tv + b \nabla f(x^*),  \, \forall  (t,b)  \in \R \times \R \}.$
The curve $\gamma(t) = x^* + t v + \varphi_a(t) \nabla f(x^*)$ defined for all $t \in \R$ such that $ta \in B(0,r)$ lies in the intersection $\mathcal{M}_{f, x^*} \cap \mathcal{L}_{f,v} $   
 and satisfies $\gamma(0)=x^*$ and 
$\gamma'(0) = v$ since $\varphi_a(0)=0$ and $\varphi'_a(0) = 0.$ By taking $v \in T_{x^*} \mathcal{M}_{f, x^*}, \, \| v \| = 1,$ and re-parametrizing the curve $\gamma$ by arc length, we get a curve satifying the hypotheses of proposition 2.

\vs 

For the submanifold $\mathcal{M}_g$ the results are obtained in the same way.  

Let $V \in \R^{n \times (n-m)}$ a matrix whose columns $v_1, \cdots, v_{n-m}$ generate the tangent space $T_{x^*} \mathcal{M}_g.$  
We define the function $\tilde{g} : \R^{n-m} \times \R^m \rightarrow \R^m$ by 
$$\tilde{g}(a,b) = g(x^* + Va +   J g(x^*)^T b),$$
for $a \in \R^{n-m}$ and $b \in \R^m.$
The Jacobian matrix $J_b \tilde{g}(0,0)$ of the function $\tilde{g}$ relatively to the variable $b$ is given by 
$$J_b \tilde{g}(0,0) = J g(x^*) \, J g(x^*)^T,$$ and thus, since $J g(x^*)$ is full rank,        $J_b \tilde{g}(0,0)$ is invertible. 
We can apply again the implicit function theorem : there exist an open ball  $B(0, \overline{r} ) \subset \R^{n-m}$ of  radius $\overline{r} > 0,$  an open ball  $B(0, \overline{r}') \subset \R^m$ of  radius $\overline{r}' > 0,$  and an implicit function $\psi : B(0,\overline{r}) \rightarrow B(0,\overline{r}')$ such that 
$$(a,b) \in B(0,\overline{r}) \times B(0,\overline{r}') \; \textrm{ and } \; \tilde{g}(a,b) = 0  \Leftrightarrow a \in B(0,\overline{r}) \; \textrm{ and } \; b= \psi(a).$$  
Again, the implicit function $\psi$ verifies $\psi(0)=0,$ and $J \psi(0)=0.$ 

Finally, the restriction of the implicit function $\psi$ to a line $\R v$ generated by a fixed vector $v \in T_{x^*} \mathcal{M}_g$ defines an implicit curve in the intersection $\mathcal{M}_g \cap \mathcal{L}_{g,v},$ where $\mathcal{L}_{g,v}$ is the affine subspace  $\mathcal{L}_{g,v} = \{ x^* + tv + J g(x^*)^T b,  \, \forall  (t,b)  \in \R \times \R^m \}.$
\vs 

\section{Second-order sufficient conditions}
In order to show the sufficient second-order sufficient conditions, we  consider the function $F : B(0, \overline{r} ) \rightarrow \R,$ defined by   
\begin{equation} F(a)=f(x^* + Va +   J g(x^*)^T \psi(a)), \label{eq71} \end{equation} where as above,   
 the implicit function $\psi : B(0,\overline{r}) \rightarrow B(0,\overline{r}'),$ with $B(0,\overline{r}) \subset \R^{n-m}$ 
and $B(0,\overline{r}') \subset \R^{m},$ verifies 
 \begin{equation} \tilde{g}(a,\psi(a)) = g(x^* + Va +   J g(x^*)^T \psi(a))=0, \, \forall a \in B(0, \overline{r} ) \label{eq8}.\end{equation}
Using the properties of $\psi$, we get $\nabla F(0)^T= \nabla f(x^*)^T V.$  But, from  first-order optimality conditions, we have $T_{x^*} \mathcal{M}_g \subset T_{x^*} \mathcal{M}_{f, x^*},$ and thus $\nabla F(0) = 0,$ since the columns of $V$ are vectors of the tangent space  $T_{x^*} \mathcal{M}_g.$   

\vs 

By derivating twice we get the Hessian matrix  $\nabla^2 F(0)$ : 
\begin{equation} \nabla^2 F(0) = V^T \nabla^2 f(x^*) V + \sum_{j=1}^m \nabla f(x^*)^T \nabla g_j(x^*) \; \nabla^2 \psi_j(0),\label{eq81}\end{equation}
where  $\nabla^2 \psi_j$ is the Hessian matrix of the $j$th component of the implicit function $\psi$ 
($\psi=(\psi_1, \cdots , \psi_m)$).  
 
In the same way, by derivating equation \eqref{eq8} twice, we get,  
 for $a=0$, \begin{equation} V^T \nabla^2 g_i(x^*) V + \sum_{j=1}^m \nabla g_i(x^*)^T \nabla g_j(x^*) \; \nabla^2 \psi_j(0) = 0, \, \textrm{ for }  i=1, \cdots, m. \label{eq82}\end{equation}    

The Hessian matrix $\nabla^2 F(0)$ is related to the Hessian matrix $\nabla_{xx}^2 L(x^*, \lambda^*).$ 
      
\begin{lemma} If $x^* \in \mathcal{M}_g$ is a local minimum and $\lambda^* \in \R^m$ the Lagrange multipliers vector associated to $x^*$,  we have the following relation : 

$$\nabla^2 F(0) = V^T \nabla_{xx}^2 L(x^*, \lambda^*) V.$$        
\end{lemma}
\begin{proof}

For any vector $a \in \R^{n-m},$ we get from \eqref{eq82} 
$$a^T V^T \nabla^2 g_i(x^*) V a  + \sum_{j=1}^m \nabla g_i(x^*)^T \nabla g_j(x^*) \; a^T \nabla^2 \psi_j(0) a = 0,  \, \textrm{ for } i=1, \cdots,m,$$ which gives, in vector form, 
$$a^T V^T \nabla^2 g(x^*) V a + J g(x^*) J g (x^*)^T \; a^T \nabla^2 \psi(0) a = 0,$$
where $a^T V^T \nabla^2 g(x^*) V a \in \R^m$ is the column vector 
with coefficients \\ $a^T V^T \nabla^2 g_i(x^*) V a, \, i=1, \cdots, m,$ and  $a^T \nabla^2 \psi(0) a \in \R^m$ the column vector 
with coefficients $a^T  \nabla^2 \psi_i(0) a, \, i=1, \cdots, m.$ From this last equality we deduce that 
\begin{equation}  a^T \nabla^2 \psi(0) a = - \, (J g(x^*) J g (x^*)^T)^{-1} 
\, a^T V^T \nabla^2 g(x^*) V a. \label{eq83} \end{equation}
But equation \eqref{eq81} leads to 
$$a^T \nabla^2 F(0) a = a^T V^T \nabla^2 f(x^*) V a + \nabla f(x^*)^T  \, J g(x^*)^T  \; a^T \nabla^2 \psi(0) a.$$ 
Using \eqref{eq83}, we obtain  

\begin{align*}
a^T \nabla^2 F(0) a  &= a^T V^T \nabla^2 f(x^*) V a 
\\  &- \nabla f(x^*)^T  \, J g(x^*)^T (J g(x^*) J g (x^*)^T)^{-1} 
\, a^T V^T \nabla^2 g(x^*) V a.
\end{align*}    
Using equality \eqref{eq55} given in the proof of proposition 3., we get $a^T \nabla^2 F(0) a = a^T V^T \nabla_{xx}^2 L(x^*, \lambda^*) V a,$ and since this equality  is verified for any vector $a \in \R^{n-m},$ we can conclude. 
\end{proof}
\vs
We now give a proof of the second-order sufficient conditions   
for  $x^* \in \mathcal{M}_g$ to be a local minimum. The classical proof of this result is obtained by a reductio ad absurdum and relies on a compactness argument (see ~\cite{LY}). The proof presented here is based on the implicit function theorem thanks to which the problem is transformed into an optimality problem without constraints.   
\vs 
  
\begin{proposition} 

If $x^* \in \mathcal{M}_g$ satisfies the first-order conditions \eqref{eq2} and \begin{equation} v^T \left(\nabla^2 f(x^*) - \sum_{i=1}^m \lambda_i^* \nabla^2 g_i(x^*) \right) v >  0, \, \forall 
v \in \textrm{Ker}(J g(x^*)), v \neq 0, \label{eq84} \end{equation} then $x^*$ is a local minimum of problem \eqref{eq1}.
\end{proposition}

\begin{proof} 

We have to define a centered open ball $B(0,R) \subset \R^n$ such that for each $h \in B(0,R)$ which verifies 
$g(x^* + h)=0,$ we have $f(x^* + h) \geq f(x^*).$

We use the orthogonal decomposition $T_{x^*} \mathcal{M}_g \oplus (T_{x^*} \mathcal{M}_g)^\perp=\R^n$ to write $h =v + w,$ with $v \in T_{x^*} \mathcal{M}_g$ and      
$w \in (T_{x^*} \mathcal{M}_g)^\perp.$ As above, let $V \in \R^{n \times (n-m)}$ a matrix whose columns $v_1, \cdots, v_{n-m}$ generate the tangent space $T_{x^*} \mathcal{M}_g.$ We can suppose without restriction that the vectors are normed and orthogonal to each other, which is equivalent to $V^T V = I_{n-m}.$      

We can write $v = V a$, for some vector $a \in \R^{n-m}$ and $w= J g(x^*)^T b,$ for some vector $b \in \R^m.$

In order to consider the implicit function $\psi : B(0,\overline{r}) \rightarrow B(0,\overline{r}')$ defined above (see equation \eqref{eq8}), we define  
$R = \min(\overline{r}, 
\nu \, \overline{r}'),$ where $\nu > 0$ verifies $\| J g(x^*)^T b \| \geq \nu \, \| b \|, \, \forall b \in \R^m$ (the existence of a constant $\nu > 0$ verifying this  inequality follows from the fact that $J g(x^*)$ is full rank). We easily verify that for $h = V a + J g(x^*)^T b$ such that $\| h \| \leq R,$ we have $\| a \| \leq \overline{r}$ and 
$\| b \| \leq \overline{r}'.$               

We can now consider the second-order Taylor series at point $0$ of the function $F : B(0, \overline{r} ) \rightarrow \R$ : 
$$F(a) = F(0) + \nabla F(0)^T a + \frac{1}{2} a^T \nabla^2 F(0) a + o(\| a\|^2)= 
F(0) + \frac{1}{2} a^T \nabla^2 F(0) a + o(\| a\|^2).$$
Using Lemma 1 we obtain $$F(a) =F(0) + \frac{1}{2}  a^T V^T \nabla_{xx}^2 L(x^*, \lambda^*) V a + o(\| a\|^2).$$ 
The second-order strict inequality 
 \eqref{eq84} implies the existence of a constant $\mu > 0$ such that       
$$ a^T V^T \nabla_{xx}^2 L(x^*, \lambda^*) V a \geq \mu \, \| a \|^2, \; \forall a \in \R^{n-m}.$$ Since $\underset{a \rightarrow 0}{\lim} \frac{o(\| a\|^2)}{\| a \|^2} = 0,$ we can take $\overline{r}$ sufficiently small such that $\| a \| < \overline{r} \Rightarrow \left\vert  \frac{o(\| a\|^2)}{\| a \|^2} \right \vert < \frac{\mu}{4}.$ For each $\| a \| < \overline{r} ,$ we have 
$$ F(a) \geq  F(0) + \| a \|^2 (\frac{\mu}{2} + \frac{o(\| a\|^2)}{\| a \|^2}) \geq F(0) + \| a \|^2 \, \frac{\mu}{4}.$$ 
In conclusion, for each  $h = V a + J g(x^*)^T b$  such that $\| h \| < R,$ we have $a \in B(0,\overline{r}),$ $b \in B(0,\overline{r}'),$ and if moreover  $g(x^* + h) =0$,  according to the property of the implicit function $\psi,$   we have $b=\psi(a).$ By using the Taylor series of $F$ we deduce 
$$F(a) = f(x^* + V a + J g(x^*)^T \psi(a)) = f(x^* + h) \geq 
f(x^*) + \| a \|^2 \, \frac{\mu}{4} \geq f(x^*).$$    
This concludes the proof.
\end{proof}

\end{document}